\documentclass[12pt]{amsart}
\usepackage{graphicx} 

\DeclareMathOperator{\lcm}{lcm}   
\newtheorem{thm}{Theorem}
\newtheorem{lem}[thm]{Lemma}

\newtheorem{claim}{Claim}
\newcommand{\Z}{\mathbb Z}
\newcommand{\R}{\mathbb R}

\begin{document}
\title{BCZ map is weakly mixing}
\author{Yitwah Cheung, Anthony Quas}
\date{\today}

\begin{abstract}
The BCZ map was introduced in \cite{BCZ} by Boca, Cobeli and Zaharescu 
as a tool to study the statistical properties of Farey sequences, whose 
relation to Riemann Hypothesis dates back to Franel and Landau (see \cite{BZ}).  

Later, J. Athreya and the first author observed that the BCZ map 
arises as a Poincare section of horocycle flow, establishing in \cite{AC} 
both ergodicity as well as zero measure-theoretic entropy.  

In this article, we prove that the BCZ map is weakly mixing, answering 
the last remaining question about the BCZ map raised in \cite{BZ}.  
The proof uses a self-similarity property of the BCZ map that derives 
from a well-known fact that horocycle flow is renormalized by the 
geodesic flow, a property already observed in \cite{AC}.  

We note that the questions of mixing and rigidity remain open.  
\end{abstract}

\maketitle

The BCZ map is a description of the induced map of the horocycle flow
on the collection of unimodular two-dimensional lattices with a horizontal vector of
length at most 1. Let $\mathcal{L}_0$ denote this collection, and more generally, let 
$\mathcal{L}_b$ denote the collection of lattices with a horizontal vector of length at most
$1-b$. 

Given a lattice in $\mathcal{L}_0$, we parameterize it by a pair 
$(s,t)$ as follows. First, $s$ is the length of the shortest horizontal vector 
(so that $(s,0)$ is a point of
the lattice). Since it is a unimodular lattice, the first horizontal row 
above the $x$-axis 
containing lattice points has $y$-coordinate $1/s$. The lattice points in this (and all) 
rows are spaced $s$ apart. We let $t$ be the lattice coordinate of the rightmost 
lattice point in the first horizontal row intersected with $\{(x,y)\colon x\le 1\}$. 
The coordinates $s$ and $t$ satisfy the inequalities $0<s,t\le 1$, $s+t>1$.

Define
$$
\Omega=\{(s,t)\colon 0<s,t\le 1;\, s+t>1\}.
$$
Thus $(s,t)\in\Omega$ corresponds to the lattice 
$$
\Lambda_{(s,t)}=\{ m(s,0)+n(t,1/s) \in\R^2\colon m,n\in \Z^2\} \in \mathcal{L}_0.
$$

The BCZ map, mentioned above, is the description of the induced map 
of the horocycle flow on $\mathcal{L}_0$ expressed in these coordinates:

$$
\Phi(s,t)=(t,u)
\text{ where $u=-s+\lfloor (1+s)/t\rfloor$t.}
$$

It is straightforward to check that $\Phi$ preserves the (normalized) restriction 
of Lebesgue  measure to $\Omega$. We denote the normalized area measure by $m$.

\begin{thm}\label{thm:main}
The BCZ map is weak-mixing.
\end{thm}

We use the following abstract lemma to show that the BCZ map does not
have an eigenfunction.
The principal fact that we exploit is that the return map of the horocycle flow to 
the collection of lattices with a horizontal vector of length at most 1 is conjugate 
(by a near-identity conjugacy) to
the return map of the horocycle flow to the collection of lattices with a horizontal vector
of length at most $1-b$.

\begin{lem}
Let $X$ be a compact metric space equipped with a Borel probability measure $m$.
Let $T$ be an ergodic measure-preserving transformation of $(X,m)$. 

Suppose that there exist $\tau>0$, a sequence of subsets $(X_k)$ of $X$, positive 
integers $(N_k)$,
and a sequence of maps $(\phi_k)$ with $\phi_k\colon X\to X_k$ 
satisfying the following properties:

\begin{enumerate}
\item\label{cond:almostfull}
$m(X_k)\to 1$;
\item 
The return map, $T_k$ of $X_k$ is conjugate to $T$ by $\phi_k$;\label{cond:conj}
\item
For all $\delta>0$, $m(\{x\colon d(\phi_k(x),x)>\delta\})\to 0$;\label{cond:almostid}
\item
$m(\{x\in X_k\colon R^{(N_k)}_{k}(x)=N_k+1\})>\tau$ for each $k$,
where for $x\in X_k$, $R^{(n)}_{k}(x)$ denotes the $n$th return time to $X_k$.\label{cond:returntime}
\end{enumerate}

Then $T$ is weak mixing.

\end{lem}

\begin{proof}
Suppose for a contradiction that $f$ has an eigenfunction with eigenvalue $\alpha\ne 1$.
By ergodicity, we may assume that $f$ has absolute value 1 almost everywhere.
Let $T_k$ denote the return map of $T$ to $X_k$.
 
We have that $f_k=f\circ\phi_k^{-1}$ is an eigenfunction for $T_k$ with the same eigenvalue.
Let $m_k$ be the normalized restriction of $m$ to $X_k$ and $\|\cdot\|_k$
be the $L^1$ norm on $X_k$ with respect to $m_k$.

We make the following claim:
\begin{equation}\label{eq:toprove}
\|f_k-f\vert_{X_k}\|_k\to 0\text{ as $k\to\infty$.}
\end{equation}

Assuming this, we complete the proof as follows:

We have 
\begin{align*}
&\|f\vert_{X_k}\circ T_k^{N_k}-\alpha^{N_k}f\vert_{X_k}\|_k\\
\le\;&\|f\vert_{X_k}\circ T_k^{N_k}-f_k\circ T_k^{N_k}\|_k
+\|f_k\circ T_k^{N_k}-\alpha^{N_k}f_k\|+\|\alpha^{N_k}f_k-\alpha^{N_k}f\vert_{X_k}\|_k\\
\le\; &2\|f-f_k\|_k\to 0
\end{align*}

On the other hand, we have that $T_k^{N_k}(x)=T^{N_k+1}x$ for a subset of $X_k$ of measure
at least $\tau$. On this set, $\left|f\vert_{X_k}\circ T_k^{N_k}(x)-\alpha^{N_k}
f\vert_{X_k}(x)\right|$ takes the value $|1-\alpha|$. 
This is the required contradiction subject to a proof of 
\eqref{eq:toprove}.

To show \eqref{eq:toprove}, let $\epsilon>0$ be arbitrary. 
Let $g$ be a continuous function such that $\|f-g\|<\epsilon$.
Let $\delta$ be such that if $d(x,y)<\delta$, then $|g(x)-g(y)|<\epsilon$.
Let $k$ be sufficiently large that $m(S_k)>1-\epsilon$, where 
$S_k=\{x\colon d(\phi_k(x),x)\le \delta\}$. Clearly $m_k(\phi_k(S_k))>1-\epsilon$
also.

We now have
\begin{align*}
&|f(\phi_k^{-1}x)-f(x)|\\
\le\;&
|f(\phi_k^{-1}x)-g(\phi_k^{-1}x)|+|g(\phi_k^{-1}x)-g(x)|+|g(x)-f(x)|.
\end{align*}
Integrating the first and last terms over $m_k$ gives at most $2\|f-g\|$ and the middle
term integrates to at most $3\epsilon$.
\end{proof}

\begin{proof}[Proof of Theorem \ref{thm:main}]
We now show that the BCZ map satisfies the conditions of the lemma.
Let $b_k$ be an arbitrary sequence of positive real numbers monotonically decreasing to 0. 
We let $\Omega_k$ be $\{(s,t)\in\Omega\colon s\le 1-b_k\}$. We see that 
$m(\Omega_k)=(1-b_k)^2$, so that condition \eqref{cond:almostfull} is satisfied.

Since the BCZ map is the return map of the horocycle flow to the set of 
lattices with a horizontal
vector of length at most 1, we see that $\Omega_k$ corresponds to the set of 
lattices with a horizontal vector of length at most $1-b_k$. 
At the level of lattices, left multiplication by 
$\begin{pmatrix}1-b_k&0\\0&(1-b_k)^{-1}\end{pmatrix}$
gives a conjugacy between the return map to lattices with a horizontal
vector of length at most
1 and the return map to lattices with a horizontal vector of length at most $1-b_k$.

Expressing this conjugacy in the $\Omega$ coordinates, one has
$$
\phi_k(s,t)=\begin{cases}
((1-b_k)s,(1-b_k)t)&\text{if $s+t>(1-b_k)^{-1}$;}\\
((1-b_k)s,(1-b_k)(js+t))&\text{otherwise},
\end{cases}
$$
where $j=\min\{n\in\mathbb N\colon (1-b_k)((j+1)s+t)>1\}$.
Notice that $d(\phi_k(s,t),(s,t))\le \sqrt2\,b_k$ for $(s,t)$ belonging to a set of measure
converging to 1, so that conditions \eqref{cond:conj} and
\eqref{cond:almostid} are satisfied.

The bulk of the work is in showing that the remaining condition, \eqref{cond:returntime},
is satisfied. To understand the return condition, it is helpful to see its geometric significance 
at the level of lattices. We have already described the link between $\Omega$ and $\Lambda_0$.
There is a similar link between $\Omega_k$ and $\Lambda_{b_k}$. 

Given $n$, we can give a precise description at the level of lattices of
the collection of those lattices in $\Omega_k$ such that the $n$th return to
$\Omega_k$ is the $(n+1)$st return to $\Omega$. Recall that the horocycle flow
uniformly decreases the slopes of all lattice points at rate 1, while keeping the $x$ coordinate
fixed. Returns to $\Lambda_0$ are precisely when there is a primitive lattice point with $y$ 
coordinate 0 and $x$ coordinate in the range $(0,1]$. Returns to $\Lambda_{b_k}$ are when 
the $x$ coordinate is in the range $(0,1-b_k]$. 

If the primitive lattice points with $x$ coordinate in the range $(0,1]$ and
positive $y$ coordinate are ordered by slope,
the horocycle flow time until the first return to $\Lambda_0$ is numerically equal 
to the slope of the first 
primitive lattice point and the flow time until the $n$th return to $\Lambda_0$ is the 
slope of the $n$th primitive lattice point. Similarly, the flow time to the $n$th return to
$\Lambda_{b_k}$ is the slope of the $n$th primitive lattice point with $x$ coordinate in the
range $(0,1-b_k]$. We summarize the conclusion of this in the following lemma.

\begin{lem}\label{lem:plus1cond}
$R^{(n)}_k(s,t)$ is $n+1$ if and only if $\lambda(s,t)$ contains exactly one primitive
lattice point with $x$ coordinate in $(1-b_k,1]$ and positive $y$ coordinate
whose slope is smaller than the slope of the $n$th
primitive lattice point with $x$ coordinate in the range $(0,1-b_k]$. 
\end{lem}

This is illustrated in Figure \ref{fig:lattice}.

\begin{figure}
\includegraphics[width=3in]{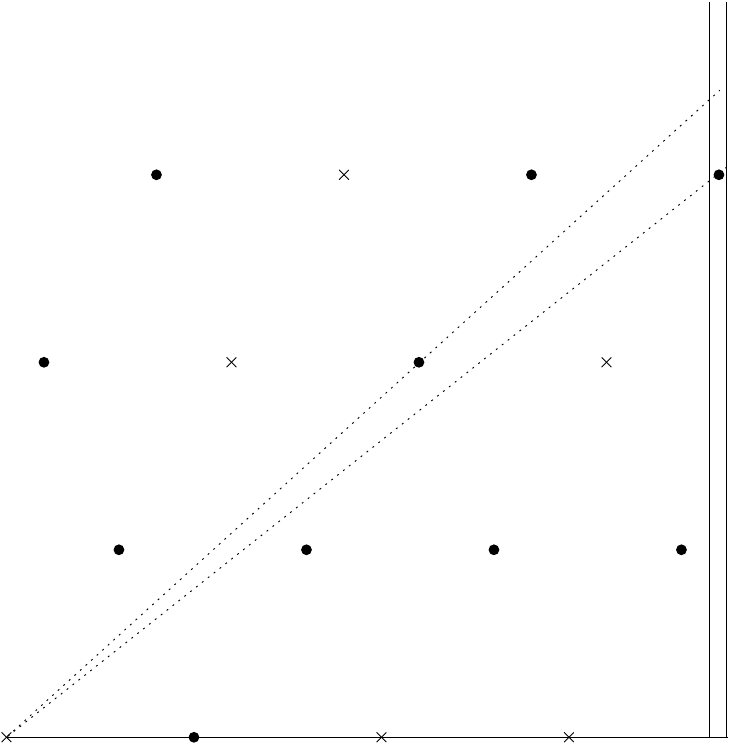}
\caption{Primitive Lattice points are marked as dots and others as crosses. 
If $x=1-b_k$ is the left
vertical line and $x=1$ is the right vertical line, then for the illustrated lattice, the 4th flow
return time to $\Lambda_0$ is the slope of the lower dotted line, while the 4th flow
return time to $\Lambda_{b_k}$ is the slope of the upper dotted line. In particular,
$R_k^{(4)}=5$ since the 4th return to $\Omega_k$ is the 5th return to $\Omega$.}\label{fig:lattice}
\end{figure}

We need to ensure that  condition \eqref{cond:returntime} is satisfied with a 
probability that is uniformly bounded away
from 0 as $k$ increases. We use the so-called second moment
method to do this. We sketch the argument.
We introduce a parameter $a$ which will later be fixed to satisfy some numerical
constraints. We then define $N_k$ to be $a/b_k$. For sufficiently large $k$, the slope of the $N_k$th primitive
point will be (from the ergodic theorem applied to the horocycle flow)
close to $\pi^2N_k/3$. We then define two boxes, each with
$x$ coordinates in the range $(1-b_k,1]$. 
The first box will have height $2a/b_k$, substantially less than $\pi^2N_k/3$, while
the second box will have height $4a/b_k$, substantially more than $\pi^2N_k/3$. We then use probabilistic
methods to show that the expected number of primitive lattice pts in $B_1$ is $\Omega(a)$ (that is,
it is bounded below, for sufficiently large, $k$ by $ca$ for a positive constant $c$ that is
independent of $k$),
while the expectation coming from points with two or more lattice points
in $B_2$ is $O(a^2)$, again where the implied constant is independent of $k$.

By choosing $a$ sufficiently small
(independent of $k$), we obtain that
with non-negligible probability there is exactly one primitive lattice point with $x$ coordinate
in $(1-b_k,1]$ below the $N_k$th lattice point with $x$ coordinate in range $(0,1-b_k]$.
This is illustrated in Figure \ref{fig:boxes}.

\begin{figure}
\includegraphics[width=3in]{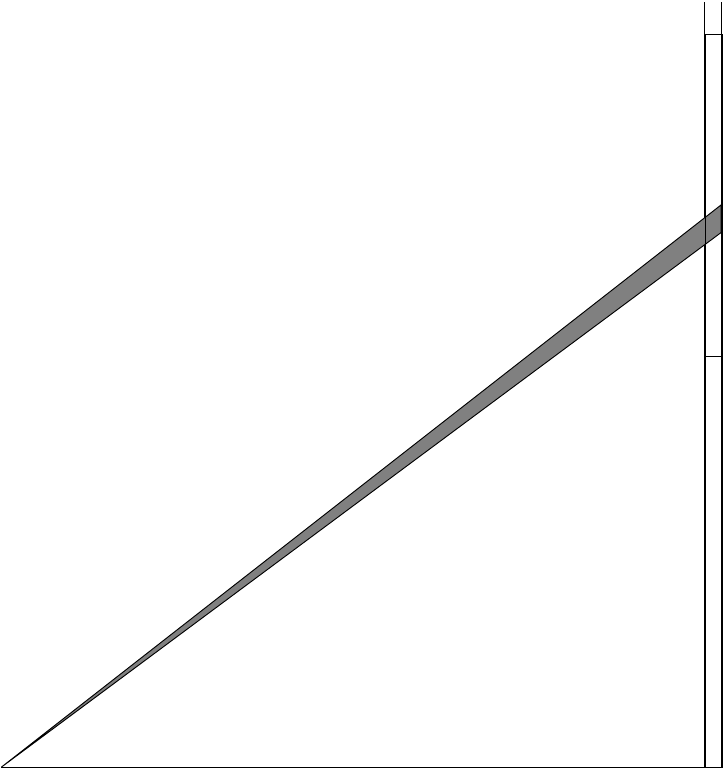}
\caption{The $N_k$th slope lies in the shaded band with overwhelming probability;
The box $B_1$ has expected number of primitive lattice points $\Omega(a)$, while
the contribution to the expected number of lattice points in $B_2$ 
coming from lattices that contain at least two lattice points is $O(a^2)$.}\label{fig:boxes}
\end{figure}

Let a positive parameter $a$ be arbitrarily chosen. We later need to pick $a$ to ensure that the
difference of the $\Omega(a)$ term and the $O(a^2)$ term remains positive. Given the parameter $a$,
we define
\begin{align*}
N_k&=\lfloor {a}/{b_k}\rfloor\\
B_{1,k}&=(1-b_k,1]\times [0,2a/b_k]\\
B_{2,k}&=(1-b_k,1]\times [0,4a/b_k].
\end{align*}

Let $F_{1,k}(s,t)$ denote the number of primitive lattice points in $\lambda(s,t)\cap B_{1,k}$.
Also, let $F_{2,k}(s,t)$ denote the number of primitive lattice points in 
$\lambda(s,t)\cap B_{2,k}$. Let $R_k=\{(s,t)\in \Omega_k\colon s\ge \frac12\}$.

We make the following claims
\begin{claim}\label{claim:one}
There exists a constant $c_1>0$ such that for all sufficiently large $k$, one has
\begin{equation*}
\int_{R_k} F_{1,k}(s,t)\,dm(s,t) \ge c_1a.
\end{equation*}
\end{claim}

\begin{claim}\label{claim:two}
There exists a constant $c_2>0$ such that for all sufficiently large $k$, one has
\begin{equation*}
\int_{R_k} F_{2,k}(s,t)\mathbf 1_{\{F_{2,k}>1\}}\,dm(s,t)
\le c_2a^2.
\end{equation*}
\end{claim}

We show how to finish the proof of the theorem assuming these claims.

Let $G_{0,k}$ be the set of  $(s,t)$ in $\Omega_{k}$ such that the $N_k$th 
horocycle flow return time of $\lambda(s,t)$ 
(to $\Lambda_{b_k}$) is between $3N_k$ and $4N_k$. 
It's well known that horocycle flow on the space of lattices is ergodic. 
The expected return time to $\Lambda_{b_k}$ is $\pi^2/(3(1-b_k)^2)$,
which lies between $3\tfrac15$ and $3\tfrac25$ for large $k$. 
As $k$ increases, so does $N_k$ and so for large $k$, the slope of the $N_k$th
primitive lattice point with $x$ coordinate in $(0,1-b_k]$ lies between $3N_k$ and $4N_k$
with probability converging to 1. That is, $m(G_{0,k})\to 1$.

We define further `good' sets as follows:

\begin{align*}
G_{1,k}&=\{(s,t)\in R_k\colon F_{1,k}(s,t)=1\} \\
G_{2,k}&=\{(s,t)\in R_k\colon F_{2,k}(s,t)=1\}.
\end{align*}
We then show that their intersection has measure bounded away from 0 and has the
properties required in \eqref{cond:returntime}.

First, we have
\begin{align*}
m(G_{1,k})&=\int_{R_k}F_{1,k}\mathbf 1_{\{F_{1,k}=1\}}\,dm\\
&=\int_{R_k}F_{1,k}-\int_{R_k}F_{1,k}\mathbf 1_{\{F_{1,k}>1\}}\,dm\\
&\ge c_1a-\int_{R_k}F_{2,k}\mathbf 1_{\{F_{2,k}>1\}}\,dm\\
&\ge c_1a - c_2a^2.
\end{align*}

Notice that $G_{1,k}\cap G_{2,k}=G_{1,k}\setminus\{(s,t)\colon F_{2,k}(s,t)>1\}$.
By Claim \ref{claim:two}, we see that the measure of the second set is at most $c_2a^2$.
In particular $G_{0,k}\cap G_{1,k}\cap G_{2,k}$ has measure bounded below by
$c_1a-3c_2a^2$ for all sufficiently large $k$. Choosing $a=c_1/(4c_2)$, we see
that $m(G_{0,k}\cap G_{1,k}\cap G_{2,k})\ge c_1^2/(16c_2)$ for all sufficiently large $k$.
By weakening the bound in the second claim if necessary, we may enlarge $c_2$
to ensure that $a<\frac1{32}$.
Finally by Lemma \ref{lem:plus1cond}, we see that if 
$(s,t)\in G_{0,k}\cap G_{1,k}\cap G_{2,k}$,
one has $R_k^{(N_k)}(s,t)=N_k+1$, establishing condition \eqref{cond:returntime} as required.

It therefore remains to verify Claims \ref{claim:one} and \ref{claim:two}.

\begin{proof}[Proof of Claim \ref{claim:one}]
For a lattice $\lambda(s,t)$, we want to count the number of primitive lattice points 
in $B_{1,k}$. The primitive lattice points of $\lambda(s,t)$ are exactly those points of
the form $-m(s,0)+n(t,1/s)$ with $\gcd(m,n)=1$. We consider the subset of 
$R_k=\{(s,t)\in\Omega_k\colon s\ge \frac12\}$
consisting of those lattices such that the $(m,n)$ lattice point lies in $B_{1,k}$.

Rewriting the constraints $-m(s,0)+n(t,\frac 1s)\in B_{1,k}$, we see 
this is the case if $(s,t)$ satisfies
\begin{align}
&1-b_k< -ms+nt\le 1\label{eq:narrow}\\
&n/s\le 2N_k.\label{eq:mn}
\end{align}

We consider only pairs $(m,n)$ with $\gcd(m,n)=1$, $N_k/2\le n\le N_k$ and $N_k/6< m\le N_k/3$.
We are also only looking at points where $s\ge 1/2$. In this range, \eqref{eq:mn} holds
automatically, so we only need to look at the $(s,t)$ pairs satisfying \eqref{eq:narrow}. 

For integers $m$ and $n$, we let $A_{m,n}=
\{(s,t)\in\Omega_k\colon s\ge \frac12\}\cap S_{m,n}$,
where 
$$
S_{m,n}=\left\{(s,t)\in \R^2\colon \frac1n -\frac{b_k}n+\frac mn s<t\le 
\frac 1n+\frac mns\right\}.
$$
That is, $A_{m,n}$ is the set of $(s,t)$ such that the $(m,n)$
lattice point of $\lambda(s,t)$ has $x$ coordinate in the range $(1-b_k,1]$:

We use a lower bound for $F_{1,k}$ consisting of a sum over $A_{m,n}$'s 
over the range described above.
$$
F_{1,k}(s,t)\ge \sum_{(m,n)\in Q}\mathbf 1_{A_{m,n}}(s,t),
$$
where $Q=\{(m,n)\in (N_k/6,N_k/3]
\times [N_k/2,N_k] \colon \gcd(m,n)=1\}$.

A typical set $A_{m,n}$ is illustrated in Figure \ref{fig:Amn}.
Notice that by the constraints on $m$ and $n$, the slope of the region is at least
$\frac 16$ and at most $\frac 23$. Since the $y$-axis intercept of the bottom line is above 0, 
a quick calculation shows that the width of the region is at least $\frac 17$, so since the 
height is $\frac {b_k}n$, we see that $m(A_{m,n})\ge 2b_k/(7n)\le b_k/(4N_k)$. 
\begin{figure}
\includegraphics[width=2in]{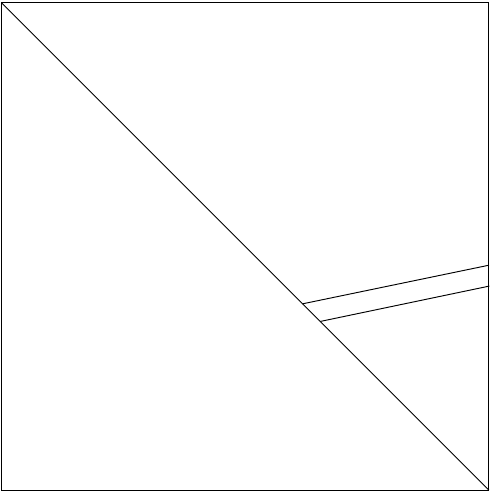}
\caption{A typical set $A_{m,n}$. The height of the strip is $b_k/n$ and the $y$ axis
intercept of the top line is $1/n$.}\label{fig:Amn}
\end{figure}
Since $(\frac{N_k}6,\frac{N_k}3]\times [\frac{N_k}2,N_k]$ contains approximately $N_k^2/12$ lattice points,
and the density of primitive points is $\frac 6{\pi^2}$, $Q$ contains approximately 
$N_k^2/(2\pi^2)$ lattice points, so at least $N_k^2/20$ for large $k$. 
Hence we obtain $\int_{R_k} F_{1,k}\,dm\ge N_kb_k/80\ge a/100$ for all large $k$.
\end{proof}

\begin{proof}[Proof of Claim \ref{claim:two}]
We first observe that $\int F_{2,k}\mathbf 1_{\{F_{2,k}>1\}}\,dm\le \int (F_{2,k}^2-F_{2,k})\,dm$
since if $k>1$ is an integer, then $k^2-k\ge k$. Since $F_{2,k}$ is a sum of indicator functions, 
$F_{2,k}^2-F_{2,k}$ can then be written as a sum of products of pairs of indicator functions (with the
diagonal terms removed).

We notice that $A_{m,n}\cap A_{m',n}=\emptyset$ for $m\ne m'$. To see this, notice that 
the intersections of $A_{m,n}$ and $A_{m',n}$
with the vertical line $x=s$ (for $s\ge \frac12$)
are intervals of height $b_k/n$ with lowest points separated by at least $1/(2n)$.

We can check that for the $(m,n)$ point of a lattice $\lambda(s,t)$
to lie in $B_2$, it is necessary that $n\le 4N_k$. We also require $m\le 2n$
as otherwise, $A_{m,n}$ does not intersect the right half of the unit square.
We make use of the simple bound
$F_{2,k}(s,t)\le \sum_{n\le 4N_k,\,m\le 2n}\mathbf 1_{A_{m,n}}(s,t)$.

Hence we see
\begin{equation}\label{eq:2bound}
\begin{split}
&\int(F_{2,k}^2-F_{2,k})\,dm\\
=&
2\sum_{n<n'\le 4N_k}\int_{R_k}
\sum_{m<2n;m'<2n'}\mathbf 1_{A_{m,n}}(s,t)\mathbf 1_{A_{m',n'}}(s,t)\,dm(s,t)\\
=&2\sum_{n<n'\le4N_k}\sum_{m<2n;m'<2n'} m(A_{m,n}\cap A_{m',n'}).
\end{split}
\end{equation}

A calculation shows that intersecting a strip with slope $s$ and height $h$
with a second strip of height $s'$ and height $h'$ gives a parallelogram
with area $hh'/|s-s'|$ and width $(h+h')/|s-s'|$.
In particular, if $n<n'\le4N_k$ then $S_{m,n}$ and $S_{m',n'}$ intersect
in a parallelogram of area $b_k^2/(mn'-m'n)$ and
width at most $b_k(1/n+1/n')/(1/nn')= b_k(n+n')\le 8N_kb_k<8a
<\frac14$.

From this, we derive a necessary condition for $A_{m,n}$ and $A_{m',n'}$ to have 
a non-trivial common intersection with $R_k$. Namely, that the lines
$t=\frac 1n+\frac mns$ and $t=\frac1{n'}+\frac{m'}{n'}s$ have an intersection
with $s$ coordinate in the range $[\frac14,\frac54]$.
In particular, taking $n<n'$, we require
$
\tfrac 14(\tfrac 1n-\tfrac 1{n'})\le \tfrac{m'}{n'}-\tfrac mn\le \tfrac 54(\tfrac 1n-\tfrac 1{n'})
$, or

\begin{equation}
\frac{4(n'-n)}5\le m'n-mn'\le 4(n'-n).\label{eq:neccond}
\end{equation}

From the first of these inequalities, we see that if there is an intersection in $R_k$, then
the difference in the slopes is equal, up to a bounded multiplicative factor, to $\frac 1n-\frac 1{n'}$.
Substituting this into the formula for the area of the parallelogram $S_{m,n}\cap S_{m',n'}$,
we see that if $A_{m,n}$ and $A_{m',n'}$ intersect, then the area of intersection is
bounded by a constant multiple of $b_k^2/(n'-n)$.

For fixed $n<n'$ we now count the number of $m\le 2n$ and $m'\le 2n'$ such that
\eqref{eq:neccond} is satisfied. Let $d=\gcd(n,n')$. Clearly for each value of $m$ and $m'$,
$m'n-mn'$ is a multiple of $d$. By the Chinese remainder theorem, the expression
$m'n-mn'$ hits each of the $nn'/d^2$ residue classes modulo $\lcm(n,n')$ 
that are multiples of $d$ exactly once as the pair  $(m,m')$ runs over $[0,n/d)\times [0,n'/d)$.
In particular, there are $(n'-n)/d$ pairs of equivalence classes of $(m,m')$
modulo $(n/d,n'/d)$ that might
satisfy \eqref{eq:neccond}. Each such class lands in the band at most $2d$ times,
as from one solution you can get to another in the same class by adding $n/d$ to $m$
and $n'/d$ to $m'$. This can be done at most $2d$ times before the constraints
$m\le 2n$ and $m'\le 2n'$ are violated.

Hence for fixed $n<n'\le 4N_k$, we have
\begin{align*}
\sum_{m\le 2n,\,m'\le 2n'}m(A_{m,n}\cap A_{m',n'})&\le \frac {cb_k^2}{n'-n}\cdot\frac{n'-n}{\gcd(n,n')}
\cdot 2\gcd(n,n')\\
&\le 2b_k^2
\end{align*}
Summing over $n$ and $n'$ in range and using \eqref{eq:2bound}, we obtain the bound
$$
\int (F_{2,k}^2-F_{2,k})\,dm\le 32cN_k^2b_k^2=: c_2a^2.
$$
\end{proof}
\end{proof}

\bibliographystyle{plain}

\end{document}